\documentclass{amsart}

\usepackage{amsmath,amsthm,amsfonts,amscd,amssymb}
\usepackage[pdftex]{graphicx}

\newtheorem{thm}{Theorem}[section]

\newtheorem{lem}[thm]{Lemma}

\theoremstyle{definition}

\theoremstyle{remark}
\newtheorem{rem}{Remark}[section]

\begin{document}

\title[Revisiting the hexagonal lattice]{Revisiting the hexagonal lattice: on optimal lattice circle packing}
\author{Lenny Fukshansky}

\address{Department of Mathematics, Claremont McKenna College, 850 Columbia Avenue, Claremont, CA 91711-6420}
\email{lenny@cmc.edu}
\subjclass{Primary: 11H06, 11H31, 52C05, 52C15, 05B40}
\keywords{lattices, circle packing, hexagonal lattice, well-rounded lattices}

\begin{abstract}
In this note we give a simple proof of the classical fact that the hexagonal lattice gives the highest density circle packing among all lattices in $\mathbb R^2$. With the benefit of hindsight, we show that the problem can be restricted to the important class of well-rounded lattices, on which the density function takes a particularly simple form. Our proof emphasizes the role of well-rounded lattices for discrete optimization problems.
\end{abstract}

\maketitle

\def\A{{\mathcal A}}
\def\B{{\mathcal B}}
\def\C{{\mathcal C}}
\def\D{{\mathcal D}}
\def\E{{\mathcal E}}
\def\F{{\mathcal F}}
\def\x{{\mathcal H}}
\def\I{{\mathcal I}}
\def\J{{\mathcal J}}
\def\K{{\mathcal K}}
\def\L{{\mathcal L}}
\def\Ll{{\mathfrak L}}
\def\M{{\mathcal M}}
\def\Mm{{\mathfrak M}}
\def\Pp{{\mathfrak P}}
\def\Aa{{\mathfrak A}}
\def\Ss{{\mathfrak S}}
\def\N{{\mathcal N}}
\def\PP{{\mathcal P}}
\def\R{{\mathcal R}}
\def\s{{\mathcal S}}
\def\V{{\mathcal V}}
\def\W{{\mathcal W}}
\def\X{{\mathcal X}}
\def\Y{{\mathcal Y}}
\def\H{{\mathcal H}}
\def\OO{{\mathcal O}}
\def\aaa{{\mathbb A}}
\def\cee{{\mathbb C}}
\def\Nn{{\mathbb N}}
\def\pee{{\mathbb P}}
\def\que{{\mathbb Q}}
\def\real{{\mathbb R}}
\def\zed{{\mathbb Z}}
\def\gmn{{\mathbb G_m^N}}
\def\qbar{{\overline{\mathbb Q}}}
\def\DL{{\underline{\Delta}}}
\def\DU{{\overline{\Delta}}}
\def\eps{{\varepsilon}}
\def\vek{{\varepsilon_k}}
\def\ahat{{\hat \alpha}}
\def\bhat{{\hat \beta}}
\def\gt{{\tilde \gamma}}
\def\h{{\tfrac12}}
\def\ba{{\boldsymbol a}}
\def\be{{\boldsymbol e}}
\def\bei{{\boldsymbol e_i}}
\def\bc{{\boldsymbol c}}
\def\bm{{\boldsymbol m}}
\def\bk{{\boldsymbol k}}
\def\bi{{\boldsymbol i}}
\def\bl{{\boldsymbol l}}
\def\bq{{\boldsymbol q}}
\def\bu{{\boldsymbol u}}
\def\bt{{\boldsymbol t}}
\def\bs{{\boldsymbol s}}
\def\bv{{\boldsymbol v}}
\def\bw{{\boldsymbol w}}
\def\bx{{\boldsymbol x}}
\def\bX{{\boldsymbol X}}
\def\bz{{\boldsymbol z}}
\def\bwy{{\boldsymbol y}}
\def\bg{{\boldsymbol g}}
\def\bY{{\boldsymbol Y}}
\def\bL{{\boldsymbol L}}
\def\baa{{\boldsymbol\alpha}}
\def\bb{{\boldsymbol\beta}}
\def\bet{{\boldsymbol\eta}}
\def\bxi{{\boldsymbol\xi}}
\def\bo{{\boldsymbol 0}}
\def\bol{{\boldsymbol 1}_L}
\def\ep{\varepsilon}
\def\p{\boldsymbol\varphi}
\def\q{\boldsymbol\psi}
\def\WR{\operatorname{WR}}
\def\rank{\operatorname{rank}}
\def\aut{\operatorname{Aut}}
\def\lcm{\operatorname{lcm}}
\def\sgn{\operatorname{sgn}}
\def\spn{\operatorname{span}}
\def\md{\operatorname{mod}}
\def\Norm{\operatorname{Norm}}
\def\dim{\operatorname{dim}}
\def\det{\operatorname{det}}
\def\Vol{\operatorname{Vol}}
\def\rk{\operatorname{rk}}
\def\md{\operatorname{mod}}
\def\sqp{\operatorname{sqp}}
\def\Aut{\operatorname{Aut}}
\def\GL{\operatorname{GL}}
\def\Sim{\operatorname{Sim}}

\section{Introduction}

The classical circle packing problem asks for an arrangement of nonoverlapping circles in $\real^2$ so that the largest possible proportion of the space is covered by them. This problem has a long and fascinating history with its origins in the works of Albrecht D\"{u}rer and Johannes Kepler. The answer to this is now known: the largest proportion of the real plane, about 90.7\%, is covered by the arrangement of circles with centers at the points of the hexagonal lattice. The first claim of a proof was made by Axel Thue in 1892, and then once again in 1910. It is generally believed however that the first complete flawless proof was produced only in 1940 by L\'{a}szl\'{o} Fejes-T\'{o}th (see \cite{conway}, \cite{rogers} for detailed accounts and bibliography). On the other hand, the fact that the hexagonal lattice gives the maximal possible circle packing density among all {\it lattice} arrangements has been known much earlier: all the necessary ingredients for the first such proof were present already in the work of Lagrange, although he himself, while aware of the circle packing problem, may not have realized that he essentially had a proof for the optimal lattice packing in hands. In fact, the notion of a lattice has not been formally introduced until the work of Gauss in 1831. A detailed history and overview of these and other developments in the direction of the circle packing problem and its (much more difficult) three-dimensional analogue, the Kepler's conjecture, can be found in the excellent recent book of G. G. Szpiro \cite{kepler}.
\smallskip

\begin{figure}[t]
\centering
\includegraphics{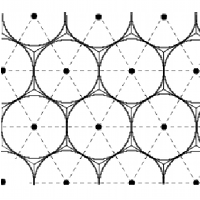}
\caption{Hexagonal lattice with Voronoi cell translates and associated circle packing}\label{fig:hex}
\end{figure}

In this note we concentrate on the lattice circle packing problem. Let us first set up the basic notation and describe the problem. Recall that a lattice $\Lambda$ in $\real^2$ is a free $\zed$-module of rank two, so $\Lambda = X\zed^2$ for some matrix $X = (\bx_1\ \bx_2) \in \GL_2(\real)$, where the column vectors $\bx_1,\bx_2$ of $X$ form a basis for $\Lambda$ and $X$ is referred to as the corresponding basis matrix. The determinant of $\Lambda$, denoted by $\det(\Lambda)$, is defined to be $|\det(X)|$, which does not depend on the particular choice of a basis for $\Lambda$. Let us now construct a circle packing associated to $\Lambda$. Define the {\it Voronoi cell} of $\Lambda$ to~be
$$\V(\Lambda) = \{ \bwy \in \real^2 : \| \bwy \| \leq \| \bwy - \bx \|\ \forall\ \bx \in \Lambda \},$$
where we write $\|\ \|$ for the Euclidean norm on $\real^2$. In other words, $\V(\Lambda)$ is the closure of the set of all vectors in the real plane which are closer to $\bo$ than to any other vector of $\Lambda$. The area of the Voronoi cell is equal to $\det(\Lambda)$, and
$$\real^2 = \bigcup_{\bwy \in \Lambda} \V(\Lambda) + \bwy,$$
meaning that the real plane is tiled with the translates of $\V(\Lambda)$. Moreover, as is clear from the definition, the interiors of these translates are disjoint. Let us inscribe a circle into each translate $\V(\Lambda) + \bwy$ of this Voronoi cell by a point of the lattice, and write $r(\Lambda)$ for the radius of this circle. No two such circles overlap, and so we have a circle packing in $\real^2$, called the {\it lattice packing} corresponding to $\Lambda$. The density of this circle packing is now given by
$$\Delta(\Lambda) = \frac{\text{area of one circle}}{\text{area of the Voronoi cell}} = \frac{\pi r(\Lambda)^2}{\det(\Lambda)}.$$
The lattice packing problem in $\real^2$ is to maximize this density function on the space of all lattices. The answer has been known since the end of the nineteenth century (see Figure \ref{fig:hex}): this density function $\Delta$ on lattices in $\real^2$ is maximized by the hexagonal lattice
$$\Lambda_h := \left( \begin{matrix} 1&\frac{1}{2} \\ 0&\frac{\sqrt{3}}{2} \end{matrix} \right) \zed^2.$$
Here we will present a proof of this fact, emphasizing the particular properties of $\Lambda_h$ that make it a solution to this optimization problem. 
\smallskip

Let us say that two lattices $\Lambda$ and $\Omega$ in $\real^2$ are {\it similar} if there exists a real constant $\alpha$ and an $2 \times 2$ orthogonal real matrix $U$ such that
$$\Omega = \alpha U \Lambda,$$
in other words, if $\Omega$ can be obtained from $\Lambda$ by rotation and dilation. Similarity is readily seen to be an equivalence relation, and it is easy to notice that the packing density function $\Delta$ is constant on each similarity class. We will prove the following classical result.

\begin{thm} \label{hexag} Let $\Lambda$ be a lattice of rank 2 in $\real^2$. Then 
\begin{equation}
\label{dense_hex}
\Delta(\Lambda) \leq \Delta(\Lambda_h) = \frac{\pi}{2\sqrt{3}} = 0.906899\dots
\end{equation}
with equality in (\ref{dense_hex}) if and only if $\Lambda$ is similar to $\Lambda_h$.
\end{thm}
\bigskip

\section{Background and standard notation}

We start by setting up some additional notation. Let $B$ be the unit circle centered at the origin in $\real^2$. Given a lattice $\Lambda$, we define Minkowski successive minima $\lambda_1 \leq \lambda_2$ of $\Lambda$ to be
$$\lambda_i = \inf \left\{ \lambda \in \real_{>0} : \Lambda \cap \lambda B \text{ contains } i \text{ linearly independent nonzero vectors} \right\},$$
where $i=1,2$. By definition of the Voronoi cell of $\Lambda$, its in-radius is equal to one half of the distance from the origin to the nearest (with respect to Euclidean norm) point of $\Lambda$, which is precisely $\lambda_1/2$, and so
\begin{equation}
\label{delta}
\Delta(\Lambda) = \frac{\pi \lambda_1^2}{4 \det(\Lambda)}.
\end{equation}
We will say that the vectors $\bx_1, \bx_2 \in \Lambda$ {\it correspond to successive minima} $\lambda_1,\lambda_2$ if they are linearly independent and
$$\|\bx_1\| = \lambda_1,\ \|\bx_2\| = \lambda_2.$$
Notice that if $\bx_1,\bx_2$ correspond to successive minima in $\Lambda$, then so do $\pm \bx_1, \pm \bx_2$. From now on, when we refer to vectors corresponding to successive minima in a lattice in $\real^2$, we will always mean a pair of such vectors so that the angle $\theta$ between them is in the interval $[0,\pi/2]$. Therefore $\cos \theta > 0$, and so
\begin{equation}
\label{xt0}
\bx_1^t \bx_2 = \| \bx_1 \| \| \bx_2 \| \cos \theta > 0.
\end{equation}
A lattice $\Lambda \subset \real^2$ is called {\it well-rounded}, abbreviated WR, if its successive minima $\lambda_1$ and $\lambda_2$ are equal. The hexagonal lattice $\Lambda_h$ is an example of a WR lattice with $\lambda_1 = \lambda_2 = 1$. Well-rounded lattices are very important in coding theory \cite{esm} and discrete optimization problems \cite{martinet}; they also come up in the context of some number theoretic problems, such as Minkowski's conjecture \cite{mcmullen} and the linear Diophantine problem of Frobenius \cite{me:sinai}. For a detailed study of the distribution of certain types of WR lattices in $\real^2$ see \cite{me:wr1} and \cite{me:wr}. In Lemma \ref{sim} below we show that the WR property is preserved under similarity, i.e. a well-rounded lattice in $\real^2$ can only be similar to another well-rounded lattice, and give a simple necessary and sufficient criterion for two WR lattices in $\real^2$ to be similar. Thus Theorem \ref{hexag} implies right away that only a WR lattice can maximize lattice packing density.
\smallskip

Our proof of Theorem \ref{hexag} emphasizes the importance of WR lattices. Specifically, we first prove that $\Delta$ must achieve its maximum at a WR lattice, hence this optimization problem can be restricted to WR lattices only. Next we show that if $\Lambda$ is WR, then $\Delta(\Lambda)$ is given by a particularly simple expression, and maximizing it becomes an easy problem. Our argument is self-contained and requires no background beyond linear algebra. For further topics in the fascinating subject of lattice packing in dimensions two and higher see \cite{conway}, \cite{hsiang}, \cite{martinet}, \cite{rogers}, and \cite{achill}. We are now ready to proceed.
\bigskip

\section{Properties of well-rounded lattices in $\real^2$}

Our goal here is to prove that the circle packing density function on the space of all lattices in $\real^2$ achieves its maximum at the hexagonal lattice. We start with a simple, but very useful lemma.

\begin{lem} \label{vec_angles} Let $\bx_1$ and $\bx_2$ be nonzero vectors in $\real^2$ so that the angle $\theta$ between them satisfies $0 < \theta < \frac{\pi}{3}$. Then
$$\| \bx_1 - \bx_2 \| < \max \{ \| \bx_1 \|, \| \bx_2 \| \}.$$
\end{lem}

\proof
Notice that $\bx_1^t \bx_2 > 0$ by (\ref{xt0}). Then, since $\theta < \frac{\pi}{3}$,
$$\frac{1}{2} < \cos \theta = \frac{\bx_1^t \bx_2}{\|\bx_1\| \|\bx_2\|},$$
and hence
\begin{eqnarray*}
\|\bx_1-\bx_2\|^2 & = & (\bx_1-\bx_2)^t (\bx_1-\bx_2) = \|\bx_1\|^2 + \|\bx_2\|^2 - 2 \bx_1^t \bx_2 \\
& < & \|\bx_1\|^2 + \|\bx_2\|^2 - \|\bx_1\| \|\bx_2\| <  \max \{ \| \bx_1 \|, \| \bx_2 \| \}^2.
\end{eqnarray*}
\endproof

Lemma \ref{vec_angles} readily implies that the angle between vectors corresponding to successive minima in a lattice cannot be~$< \pi/3$.

\begin{lem} \label{angle} Let $\Lambda \subset \real^2$ be a lattice of full rank with successive minima $\lambda_1 \leq \lambda_2$, and let $\bx_1,\bx_2$ be the vectors in $\Lambda$ corresponding to $\lambda_1,\lambda_2$, respectively. Let $\theta \in [0,\pi/2]$ be the angle between $\bx_1$ and $\bx_2$. Then
$$\pi/3 \leq \theta \leq \pi/2.$$
\end{lem}

\proof
Assume that $\theta < \pi/3$, then Lemma \ref{vec_angles} implies that
$$\| \bx_1 - \bx_2 \| < \| \bx_2 \| = \lambda_2,$$
which contradicts the definition of $\lambda_2$ since the vectors $\bx_1$ and $\bx_1-\bx_2$ are linearly independent.
\endproof

We can now prove that vectors corresponding to successive minima in a lattice in $\real^2$ form a basis.

\begin{lem} \label{mink} Let $\Lambda$ be a lattice in $\real^2$ with successive minima $\lambda_1 \leq \lambda_2$and let $\bx_1,\bx_2$ be the vectors in $\Lambda$ corresponding to $\lambda_1,\lambda_2$, respectively. Then $\bx_1,\bx_2$ form a basis for $\Lambda$.
\end{lem}

\proof
Let $\bwy_1 \in \Lambda$ be a shortest vector extendable to a basis in $\Lambda$, and let $\bwy_2 \in \Lambda$ be a shortest vector such that $\bwy_1,\bwy_2$ is a basis of $\Lambda$. By picking $\pm \bwy_1, \pm \bwy_2$ if necessary we can ensure that the angle between these vectors is no greater than~$\pi/2$. Then 
$$0 < \|\bwy_1\| \leq \|\bwy_2\|,$$
and for any vector $\bz \in \Lambda$ with $\|\bz\| < \|\bwy_2\|$ the pair $\bwy_1,\bz$ is {\it not} a basis for $\Lambda$. Since $\bx_1,\bx_2 \in \Lambda$, there must exist integers $a_1,a_2,b_1,b_2$ such that
\begin{equation}
\label{x_y}
\left( \bx_1\ \bx_2 \right) = \left( \bwy_1\ \bwy_2 \right) \left( \begin{matrix} a_1&b_1 \\ a_2&b_2 \end{matrix} \right).
\end{equation}
Let $\theta_x$ be the angle between $\bx_1,\bx_2$, and $\theta_y$ be the angle between $\bwy_1,\bwy_2$, then $\pi/3 \leq \theta_x \leq \pi/2$ by Lemma \ref{angle}. Moreover, $\pi/3 \leq \theta_y \leq \pi/2$: indeed, suppose $\theta_y < \pi/3$, then by Lemma \ref{vec_angles},
$$\|\bwy_1-\bwy_2\| < \|\bwy_2\|,$$
however $\bwy_1, \bwy_1-\bwy_2$ is a basis for $\Lambda$ since $\bwy_1, \bwy_2$ is; this contradicts the choice of~$\bwy_2$. Define 
$$\D = \left| \det \left( \begin{matrix} a_1&b_1 \\ a_2&b_2 \end{matrix} \right) \right|,$$
then $\D$ is a positive integer, and taking determinants of both sides of (\ref{x_y}), we obtain
\begin{equation}
\label{x_y1}
\|\bx_1\| \|\bx_2\| \sin \theta_x = \D \|\bwy_1\| \|\bwy_2\| \sin \theta_y.
\end{equation}
Notice that by definition of successive minima, $\|\bx_1\| \|\bx_2\| \leq \|\bwy_1\| \|\bwy_2\|$, and hence (\ref{x_y1}) implies that
$$\D = \frac{\|\bx_1\| \|\bx_2\|}{\|\bwy_1\| \|\bwy_2\|} \frac{\sin \theta_x}{\sin \theta_y} \leq \frac{2}{\sqrt{3}} < 2,$$
meaning that $\D = 1$. Combining this observation with (\ref{x_y}), we see that
$$\left( \bx_1\ \bx_2 \right) \left( \begin{matrix} a_1&b_1 \\ a_2&b_2 \end{matrix} \right)^{-1} = \left( \bwy_1\ \bwy_2 \right),$$
where the matrix $\left( \begin{matrix} a_1&b_1 \\ a_2&b_2 \end{matrix} \right)^{-1}$ has integer entries. Therefore $\bx_1,\bx_2$ is also a basis for $\Lambda$, completing the proof.
\endproof

\begin{rem} We note that if we replace $\real^2$ with $\real^d$ then the statement of Lemma \ref{mink} is no longer true for $d \geq 5$ (see for instance \cite{pohst}).
\end{rem}

We will call a basis for a lattice as in Lemma \ref{mink} a {\it minimal basis}. The goal of the next three lemmas is to show that the lattice packing density function $\Delta$ attains its maximum in $\real^2$ on the set of well-rounded lattices.

\begin{lem} \label{WR_reduction} Let $\Lambda$ and $\Omega$ be lattices of full rank in $\real^2$ with successive minima $\lambda_1(\Lambda),\lambda_2(\Lambda)$ and $\lambda_1(\Omega),\lambda_2(\Omega)$ respectively. Let $\bx_1,\bx_2$ and $\bwy_1,\bwy_2$ be vectors in $\Lambda$ and $\Omega$, respectively, corresponding to successive minima. Suppose that $\bx_1=\bwy_1$, and angles between the vectors $\bx_1,\bx_2$ and $\bwy_1,\bwy_2$ are equal, call this common value $\theta$. Suppose also that
$$\lambda_1(\Lambda) = \lambda_2(\Lambda).$$
Then
$$\Delta(\Lambda) \geq \Delta(\Omega).$$
\end{lem}

\proof
By Lemma \ref{mink}, $\bx_1,\bx_2$ and $\bwy_1,\bwy_2$ are minimal bases for $\Lambda$ and $\Omega$, respectively. Notice that
\begin{eqnarray*}
\lambda_1(\Lambda) & = & \lambda_2(\Lambda) = \|\bx_1\| = \|\bx_2\| \\
& = & \|\bwy_1\| = \lambda_1(\Omega)\leq \|\bwy_2\| = \lambda_2(\Omega).
\end{eqnarray*}
Then, by (\ref{delta}),
\begin{eqnarray}
\label{WR_pack}
\Delta(\Lambda) & = & \frac{\pi \lambda_1(\Lambda)^2}{4 \det(\Lambda)} = \frac{ \lambda_1(\Lambda)^2 \pi}{4 \|\bx_1\| \|\bx_2\| \sin \theta} = \frac{ \pi}{4 \sin \theta} \nonumber \\
& \geq & \frac{ \lambda_1(\Omega)^2 \pi}{4 \|\bwy_1\| \|\bwy_2\| \sin \theta} = \frac{ \lambda_1(\Omega)^2 \pi}{4 \det(\Omega)} = \Delta(\Omega).
\end{eqnarray}
\endproof

The following lemma is a converse to Lemma \ref{angle}.

\begin{lem} \label{angle1} Let $\Lambda \subset \real^2$ be a lattice of full rank, and let $\bx_1, \bx_2$ be a basis for $\Lambda$ such that
$$\|\bx_1\| = \|\bx_2\|,$$
and the angle $\theta$ between these vectors lies in the interval $[\pi/3, \pi/2]$. Then $\bx_1,\bx_2$ is a minimal basis for $\Lambda$. In particular, this implies that $\Lambda$ is WR.
\end{lem}

\proof
Let $\bz \in \Lambda$, then $\bz = a\bx_1+b\bx_2$ for some $a,b \in \zed$. Then
$$\|\bz\|^2 =  a^2 \|\bx_1\|^2 + b^2 \|\bx_2\|^2 + 2ab \bx_1^t \bx_2 = (a^2+b^2 + 2ab \cos \theta) \|\bx_1\|^2.$$
If $ab \geq 0$, then clearly $\|\bz\|^2 \geq \|\bx_1\|^2$. Now suppose $ab < 0$, then again
$$\|\bz\|^2 \geq (a^2+b^2 - |ab|) \|\bx_1\|^2 \geq \|\bx_1\|^2,$$
since $\cos \theta \leq 1/2$. Therefore $\bx_1,\bx_2$ are shortest nonzero vectors in $\Lambda$, hence they correspond to successive minima, and so form a minimal basis. Thus $\Lambda$ is WR, and this completes the proof.
\endproof

\begin{lem} \label{WR_const} Let $\Lambda$ be a lattice in $\real^2$ with successive minima $\lambda_1,\lambda_2$ and corresponding basis vectors $\bx_1,\bx_2$, respectively. Then the lattice
$$\Lambda_{\WR} = \left( \bx_1\ \frac{\lambda_1}{\lambda_2} \bx_2 \right) \zed^2$$
is WR with successive minima equal to $\lambda_1$.
\end{lem}

\proof
By Lemma \ref{angle}, the angle $\theta$ between $\bx_1$ and $\bx_2$ is in the interval $[\pi/3,\pi/2]$, and clearly this is the same as the angle between the vectors $\bx_1$ and $\frac{\lambda_1}{\lambda_2} \bx_2$. Then by Lemma \ref{angle1}, $\Lambda_{\WR}$ is WR with successive minima equal to $\lambda_1$.
\endproof

Now combining Lemma \ref{WR_reduction} with Lemma \ref{WR_const} implies that 
\begin{equation}
\label{WR_delta}
\Delta(\Lambda_{\WR}) \geq \Delta(\Lambda)
\end{equation}
for any lattice $\Lambda \subset \real^2$, and (\ref{WR_pack}) readily implies that the equality in (\ref{WR_delta}) occurs if and only if $\Lambda = \Lambda_{\WR}$, which happens if and only if $\Lambda$ is well-rounded. Therefore the maximum packing density among lattices in $\real^2$ must occur at a WR lattice, and so for the rest of this section we talk about WR lattices only. Next observation is that for any WR lattice $\Lambda$ in $\real^2$, (\ref{WR_pack}) implies:
$$\sin \theta = \frac{\pi}{4 \Delta(\Lambda)},$$
meaning that $\sin \theta$ is an invariant of $\Lambda$, and does not depend on the specific choice of the minimal basis. Since by our conventional choice of the minimal basis and Lemma \ref{angle}, this angle $\theta$ is in the interval $[ \pi/3, \pi/2]$, it is also an invariant of the lattice, and we call it the {\it angle of $\Lambda$}, denoted by $\theta(\Lambda)$. 

\begin{lem} \label{sim} Let $\Lambda$ be a WR lattice in $\real^2$. A lattice $\Omega \subset \real^2$ is similar to $\Lambda$ if and only if $\Omega$ is also WR and $\theta(\Lambda) = \theta(\Omega).$
\end{lem}

\proof
First suppose that $\Lambda$ and $\Omega$ are similar. Let $\bx_1,\bx_2$ be the minimal basis for $\Lambda$. There exist a real constant $\alpha$ and a real orthogonal $2 \times 2$ matrix $U$ such that $\Omega = \alpha U \Lambda$. Let $\bwy_1,\bwy_2$ be a basis for $\Omega$ such that
$$(\bwy_1\ \bwy_2) = \alpha U (\bx_1\ \bx_2).$$
Then $\|\bwy_1\| = \|\bwy_2\|$, and the angle between $\bwy_1$ and $\bwy_2$ is $\theta(\Lambda) \in [\pi/3, \pi/2]$. By Lemma \ref{angle1} it follows that $\bwy_1,\bwy_2$ is a minimal basis for $\Omega$, and so $\Omega$ is WR and $\theta(\Omega) = \theta(\Lambda)$.

Next assume that $\Omega$ is WR and $\theta(\Omega) = \theta(\Lambda)$. Let $\lambda(\Lambda)$ and $\lambda(\Omega)$ be the respective values of successive minima of $\Lambda$ and $\Omega$. Let $\bx_1,\bx_2$ and $\bwy_1,\bwy_2$ be the minimal bases for $\Lambda$ and $\Omega$, respectively. Define
$$\bz_1 = \frac{\lambda(\Lambda)}{\lambda(\Omega)} \bwy_1,\ \bz_2 = \frac{\lambda(\Lambda)}{\lambda(\Omega)} \bwy_2.$$
Then $\bx_1,\bx_2$ and $\bz_1,\bz_2$ are pairs of points on the circle of radius $\lambda(\Lambda)$ centered at the origin in $\real^2$ with equal angles between them. Therefore, there exists a $2 \times 2$ real orthogonal matrix $U$ such that
$$(\bwy_1\ \bwy_2) = \frac{\lambda(\Lambda)}{\lambda(\Omega)} (\bz_1\ \bz_2) = \frac{\lambda(\Lambda)}{\lambda(\Omega)}  U (\bx_1\ \bx_2),$$
and so $\Lambda$ and $\Omega$ are similar lattices. This completes the proof.
\endproof

We are now ready to prove the main result.

\proof[Proof of Theorem \ref{hexag}]
The density inequality (\ref{WR_delta}) says that the largest lattice packing density in $\real^2$ is achieved by some WR lattice $\Lambda$, and (\ref{WR_pack}) implies that 
\begin{equation}
\label{hx1}
\Delta(\Lambda) = \frac{\pi}{4 \sin \theta(\Lambda)},
\end{equation}
meaning that a smaller $\sin \theta(\Lambda)$ corresponds to a larger $\Delta(\Lambda)$. Lemma \ref{angle} implies that $\theta(\Lambda) \geq \pi/3$, meaning that $\sin \theta(\Lambda) \geq \sqrt{3}/2$. Notice that if $\Lambda$ is the hexagonal lattice
$$\Lambda_h = \left( \begin{matrix} 1&\frac{1}{2} \\ 0&\frac{\sqrt{3}}{2} \end{matrix} \right) \zed^2,$$
then $\sin \theta(\Lambda) = \sqrt{3}/2$, meaning that the angle between the basis vectors $(1,0)$ and $(1/2, \sqrt{3}/2)$ is $\theta = \pi/3$, and so by Lemma \ref{angle1} this is a minimal basis and $\theta(\Lambda) = \pi/3$. Hence the largest lattice packing density in $\real^2$ is achieved by the hexagonal lattice. This value now follows from (\ref{hx1}).

Now suppose that for some lattice $\Lambda$, $\Delta(\Lambda) = \Delta(\Lambda_h)$, then by (\ref{WR_delta}) and a short argument after it $\Lambda$ must be WR, and so
$$\Delta(\Lambda) = \frac{\pi}{4 \sin \theta(\Lambda)} = \Delta(\Lambda_h) =  \frac{\pi}{4 \sin \pi/3}.$$
Then $\theta(\Lambda) = \pi/3$, and so $\Lambda$ is similar to $\Lambda_h$ by Lemma \ref{sim}. This completes the proof.
\endproof
\bigskip

{\bf Acknowledgment.} I would like to thank Sinai Robins for his helpful comments on the subject of this paper.
\bigskip

\bibliographystyle{plain}  
\bibliography{fukshansky_paper}    

\end{document}